\documentclass[preprint,10pt]{elsarticle}

\usepackage{amsmath, amssymb}
\usepackage{graphicx}
\journal{Acta Mathematica Scientia}

\begin{document}

\begin{frontmatter}

\title{Solvability of a non-local problem with integral transmitting condition for mixed type equation with Caputo fractional derivative}

\author{P.Agarwal}

\address{Department of Mathematics,
Anand International College of Engineering,
Jaipur, India}
\ead{goyal.praveen2011@gmail.com}

\author{A.S.Berdyshev}

\address{Kazakh National Pedagogical University named after Abai, Almaty, Kazakhstan}
\ead{berdyshev@mail.ru}

\author{E.T.Karimov}
\address{Institute of Mathematics, National University of Uzbekistan, Tashkent, Uzbekistan}
\ead{erkinjon@gmail.com}

\begin{abstract}
In the present paper, we discuss solvability questions of a non-local problem with integral form transmitting conditions for diffusion-wave equation with the Caputo fractional derivative in a domain bounded by smooth curves. The uniqueness of the solution of the formulated problem we prove using energy integral method with some modifications. The existence of solution will be proved by equivalent reduction of the studied problem into a system of second kind Fredholm integral equations.
\end{abstract}


\end{frontmatter}

\section*{1. Introduction}

It is well-known that partial differential equations (PDEs) are on the base of many mathematical models for real-life processes. On an application of mixed type PDEs for the first time was mentioned by S.A.Chaplygin [1]. Later, fundamental results on this direction were obtained by F.Tricomi, S.Gellerstedt, I.M.Frankl, M.A.Lavrent'ev, A.V.Bitsadze, M.H.Protter, C.Z.Morawetz and others. Theory of boundary-value problems for various mixed type PDEs is one of the continuously and intensively developing theories of modern mathematics. Omitting huge amount works, we would like just mention some works, where application of mixed parabolic-hyperbolic type equations were object of investigations. Precisely, in the work [2], a hyperbolic-parabolic system arising in pulse combustion is investigated. Work by Ostertagova et al. [3] is devoted to the mechanical problems, described by second order parabolic-hyperbolic equations.

Regarding the studying of various type boundary-value problems for parabolic-hyperbolic type equations we refer the readers to the monograph [4]. 

Regarding the investigations on parabolic-hyperbolic equations with three lines of type-changing we note works [5-7]. 

We as well would like to note results on local and non-local problems for parabolic-hyperbolic type equations with fractional order derivatives. Precisely, in [8] the Tricomi and Gellerstedt problem for parabolic-hyperbolic equation with the Riemann-Liouvill fractional operator in the hyperbolic part were under discussion and unique solvability of these problems were proved. In [9] authors consider the same equation, but with two lines of type-changing in a domain with deviation from the characteristics. In above-mentioned works, authors used special transmitting conditions on type-changing lines. In [10, 11], transmitting conditions were generalized and as well instead of the Riemann-Liouville fractional differential operator, Caputo fractional differential operator was considered. Under certain assumptions on given functions and parameters, unique solvability of considered problems were proved.

In [12-14] various linear and semilinear parabolic-hyperbolic type equations were investigated by numerical methods.

In the present paper we investigate non-local problem with integral form transmitting condition for parabolic-hyperbolic equation with Caputo fractional derivative on time varibale in a domain, bounded by smooth curves. The difficulties of all three aspects: influence of fractional derivative, form of considered domain and transmitting conditions together made evaluations complicative and we have to deal with different cases, requiring careful intention of us. 

We divided paper into seven sections, giving in an introduction brief review of papers, related to the topic of this paper. In next sections we formulate problem and obtain main functional correlations. The uniqueness and the existence of the solution of considered problem is proved in separate sections. At the end we give conclusions. 

\section*{2. Preliminaries}

In this section we give some known facts, which we need further.

\subsection*{2.1. The Cauchy problem for wave equation:}
A function $u(\xi,\eta)\in C(\overline{\Omega})\cap C^2(\Omega)$ we call as a solution of the Cauchy problem for the wave equation, if it satisfies equation
$$
\frac{\partial ^2 u(\xi,\eta)}{\partial \xi \partial \eta}=f(\xi,\eta)
$$
and initial conditions
$$
u(\xi,\eta)|_{\xi=\eta}=\tau(\eta),\,\,\,\left(\frac{\partial u(\xi,\eta)}{\partial \xi}-\frac{\partial u(\xi,\eta)}{\partial \eta}\right)|_{\xi=\eta}=\nu(\eta).
$$
This solution can be represented as [15]
$$
u(\xi,\eta)=\tau(\xi)+\tau(\eta)-\int\limits_{\xi}^{\eta}\nu(z)dz-\int\limits_{\xi}^{\eta} d\xi_1\int\limits_{\xi_1}^{\eta} f(\xi_1,\eta_1)d\eta_1. \eqno (1)
$$
\subsection*{2.2. First boundary problem (FBP) for heat equation with Caputo derivative:}
We consider the following heat equation with the Caputo fractional derivative:
$$
\frac{\partial ^2 u(x,t)}{\partial x^2}-_CD_{0t}^\lambda u(x,t)=f(x,t), \eqno (2)
$$
where $0<\lambda \leq 1$,
$$
_CD_{0t}^\lambda g(t)=\left\{
\begin{array}{l}
\frac{1}{\Gamma (1-\lambda)}\int\limits_0^t \frac{g'(z)}{(t-z)^\lambda}dz, \,\,\,\, \emph{if} \,\,\,\,0<\lambda<1,\\
\frac{d g(t)}{dt},\,\,\,\,\,\,\,\,\,\,\,\, \,\,\,\,\,\,\,\,\,\,\,\,\,\,\,\,\,\,\,\,\,\, \,\,\,\,\,\emph{if} \,\,\,\,\,\,\,\,\,\,\,\,\,\,\,\lambda=1\\
\end{array}
\right.
$$
is the Caputo fractional differential operator of the order $\lambda$ [16].

\textbf{FBP}. To find a solution of the equation (2) in a rectangular domain\\
 $\left\{(x,t):\, 0\leq x \leq 1,\,0\leq t\leq 1\right\}$, satisfying the following boundary conditions:
$$
u(x,0)=\varphi_1(x),\,\,0\leq x \leq 1,\,\,\,\,u(0,t)=\varphi_2(t),\,\,u(1,t)=\varphi_3(t),\,0\leq t \leq 1.
$$
This solution can be represented as [17]
$$
\begin{array}{l}
u(x,t)=\int\limits_0^t G_\xi(x,t,0,t_1)\varphi_2(t_1)dt_1-\int\limits_0^t G_\xi(x,t,1,t_1)\varphi_3(t_1)dt_1+\\
+\int\limits_0^1 \overline{G}(x-x_1,t)\varphi_1(x_1)dx_1-\int\limits_0^1 \int\limits_0^t G(x,t,x_1,t_1) f(x_1, t_1)dx_1 dt_1,\\
\end{array}
\eqno (3)
$$
where
$$
\overline{G}(x-x_1,t)=\frac{1}{\Gamma (1-\lambda)} \int\limits_0^t t_1^{-\lambda} G(x,t,x_1,t_1) dt_1,
$$
$$
G(x,t,x_1,t_1)=\frac{(t-t_1)^{\beta-1}}{2}\sum_{n=-\infty}^{+\infty} \left[e_{1,\beta}^{1,\beta} \left(-\frac{|x-x_1+2n|}{(t-t_1)^\beta}\right)
-e_{1,\beta}^{1,\beta} \left(-\frac{|x+x_1+2n|}{(t-t_1)^\beta}\right)\right]\eqno (4)
$$
is the Green's function of FBP for the equation (2), $\beta=\lambda/2$,
$$
e_{1,\beta}^{1,\beta} \left(z\right)=\sum_{n=0}^{+\infty} \frac{z^n}{n!\Gamma(\beta-\beta n)}
$$
is the Wright type function [17].

\subsection*{2.3. Properties of the Wright type function:}
\begin{itemize}
\item
At $\alpha>\beta, \alpha>0$ for any $z\in \mathbb{C}$ the following expression [17]
$$
e_{\alpha,\beta}^{\mu,\delta} \left(z\right)=\sum_{n=0}^{+\infty} \frac{z^n}{\Gamma(\alpha n+\mu)\Gamma(\delta-\beta n)}\eqno (5)
$$
is valid.

\item
For any $z\in \mathbb{C}$ the following relations [17]
$$
\frac{1}{\alpha}e_{\alpha,\beta}^{\mu-1,\delta}(z)+\frac{1}{\beta}e_{\alpha,\beta}^{\mu,\delta-1}(z)=\left[\frac{\mu-1}{\alpha}-\frac{\delta-1}{\beta}
\right]e_{\alpha,\beta}^{\mu,\delta}(z)\eqno (6)
$$
$$
\frac{1}{z}e_{\alpha,\beta}^{-k,\delta}(z)=e_{\alpha,\beta}^{\alpha-k,\delta-\beta}(z),\,\,
\frac{1}{z}e_{\alpha,\beta}^{\mu,-k}(z)=e_{\alpha,\beta}^{\mu+\alpha,-k-\beta}(z) \eqno (7)
$$
and the following formula of differentiation [17]
$$
\frac{d}{dz}e_{\alpha,\beta}^{\mu,\delta}(z)=-\frac{1}{\beta z}\left[e_{\alpha,\beta}^{\mu,\delta-1}(z)+(1-\delta)e_{\alpha,\beta}^{\mu,\delta}(z)\right]\eqno (8)
$$
are true.

\end{itemize}

\section*{3. Formulation of the problem }
Let $\Omega \subset R^2$ be a finite simple-connected domain (see Figure 1) and $\Omega=\Omega_0\cup \Omega_i\cup AB\cup BC\cup AD$.

We consider an equation
\[
Lu=f\left( x,t \right),\eqno (9)
\]
where
\[
Lu\equiv \left\{ \begin{array}{l}
   u_{xx}-_CD_{0t}^\lambda u,\,\,\,(x,t)\in \Omega_0,  \\
   u_{xx}-u_{tt},\,\,\,(x,t)\in \Omega_i \, (i=\overline{1,3}).  \\
\end{array} \right.
\]

Smooth curves $\gamma_1:\, t=-\gamma_1(x), \gamma_1(0)=\gamma_1(1)=0$,
$\gamma_2:\, x=-\gamma_2(t), \gamma_2(0)=\gamma_2(1)=0$,
$\gamma_3:\, x=-\gamma_3(t)+1, \gamma_3(0)=\gamma_3(1)=0$ lie strictly inside of characteristic triangles. Moreover, $\gamma_i(s)$ are twice differentiable and $s\pm \gamma_i(s)$ \, $(0\leq s\leq 1,\, i=\overline{1,3})$ are monotonically increase.

\begin{figure}
\begin{center}
  \includegraphics[width=0.7\textwidth]{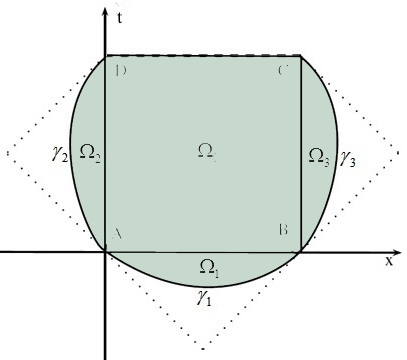}\\
  \caption{Domain $\Omega$}\label{fig1}
\end{center}
\end{figure}

We formulate the following non-local problem for the equation (9):

\textbf{Problem S.} Find a solution of the equation (9) from the following class of functions
$$
W=\left\{u:\, u(x,t)\in C(\overline{\Omega}), u_{xx}, _CD_{0t}^\lambda u \in C(\Omega_0), u(x,t)\in C^2(\Omega_i), i=\overline{1,3}\right\},
$$
first derivatives of which are continuous up to boundaries of $\Omega_0$ and $\Omega_i (i=\overline{1,3})$, satisfies non-local conditions
$$
\left[u_x-u_t\right]\left(\theta_1(s)\right)=\sigma_1\left[u_x+u_t\right]\left(\theta_1^{*}(s)\right),\,\,0<s<1, \eqno (10)
$$
$$
\left[u_x-u_t\right]\left(\theta_2(s)\right)=\sigma_2\left[u_x+u_t\right]\left(\theta_2^{*}(s)\right),\,\,0<s<1, \eqno (11)
$$
$$
\left[u_x+u_t\right]\left(\theta_3(s)\right)=\sigma_3\left[u_x-u_t\right]\left(\theta_3^{*}(s)\right),\,\,0<s<1, \eqno (12)
$$
condition on two point
$$
u(A)=u(B)=0\eqno (13)
$$
and transmitting conditions
$$
\lim\limits_{t\rightarrow+0} {_CD_{0t}^\lambda} u(x,t)=\alpha_1 u_t(x,-0)+\beta_1\int\limits_x^1u_t(s,-0) P_1(x,s)ds,\eqno (14)
$$
$$
u_x(+0,t)=\alpha_2 u_x(-0,t)+\beta_2\int\limits_t^1 u_x(-0,s) P_2(t,s)ds,\eqno (15)
$$
$$
u_x(1-0,t)=\alpha_3 u_x(1+0,t)+\beta_3\int\limits_t^1 u_x(1+0,s) P_3(t,s)ds. \eqno (16)
$$
Here $\sigma_i,\, \alpha_i,\, \beta_i$ are given real numbers such that $\alpha_i^2+\beta_i^2\neq 0$, $P_i(\cdot,\cdot)$ are given functions, $\theta_1(s),\theta_2(s),\theta_3(s) \left[\theta_1^{*}(s),\theta_2^{*}(s),\theta_3^{*}(s)\right]$ are affixes of points of intersection of curves $\gamma_i(s) (i=\overline{1,3})$ with characteristics $x-t=s, t-x=s, x+t=1+s$ $\left[x+t=s, x+t=s, t-x=1+s\right]$ of the equation (9), respectively.

\section*{4. Main functional correlations}

Based on representation of the solution of the Cauchy problem for the wave equation we can represent solution of problem S in domains $\Omega_i (i=\overline{1,3})$ as follows:
$$
u(\xi,\eta)=\frac{1}{2}\left[\tau_1^-(\xi)+\tau_1^-(\eta)-\int\limits_{\xi}^{\eta}\nu_1^-(z)dz\right]-\int\limits_{\xi}^{\eta} d\xi_1\int\limits_{\xi_1}^{\eta} f_1(\xi_1,\eta_1)d\eta_1, \eqno (17)
$$
$$
u(\xi,\eta)=\frac{1}{2}\left[\tau_2^-(\xi)+\tau_2^-(-\eta)-\int\limits_{\xi}^{-\eta}\nu_2^-(z)dz\right]-\int\limits_{\xi}^{-\eta} d\xi_1\int\limits_{\xi_1}^{-\eta} f_1(\xi_1,\eta_1)d\eta_1, \eqno (18)
$$
$$
u(\xi,\eta)=\frac{1}{2}\left[\tau_3^+(\xi-1)+\tau_3^+(1-\eta)-\int\limits_{\xi-1}^{1-\eta}\nu_3^+(z)dz\right]-\int\limits_{\xi-1}^{1-\eta} d\xi_1\int\limits_{\xi_1}^{1-\eta} f_1(\xi_1,\eta_1)d\eta_1, \eqno (19)
$$
where
$$
\begin{array}{l}
\xi=x+t,\,\eta=x-t,\,4f_1\left(\xi,\eta\right)=f\left(\frac{\xi+\eta}{2},\frac{\xi-\eta}{2}\right),\, \tau_1^\pm(x)=u(x,\pm 0),\\
\tau_2^\pm (t)=u(\pm 0,t),\, \tau_3^\pm (t)=u(1\pm 0,t),\, \nu_1^+(x)=\lim\limits_{t\rightarrow +0}{_CD_{0t}^\lambda} u(x,t),\\
\nu_1^-(x)=u_t(x,-0),\, \nu_2^\pm (t)=u_x(\pm 0,t),\, \nu_3^\pm (t)=u_x(1\pm 0,t).
\end{array}
\eqno (20)
$$
Based on conditions imposed to $\gamma_i (i=\overline{1,3})$, equations of these curves can be represented as $\xi=\rho(\eta)$ and $\eta=v(\xi)$ such that $\rho(v(\xi))=\xi$. Since $\theta_1 (\rho(s),s),\,\theta_2(\rho(-s),s),\, \theta_3(1+s,v(1+s))$ $\theta_1^{*}(s,v(s)),\,\theta_2^{*}(s,v(s)),\,\theta_3^{*}(\rho(1-s),1-s)$, then from (10) and (17), (11) and (18), (12) and (19), we obtain the following functional correlations on the lines of type-changing, upbringing from the hyperbolic parts of of the mixed domain:
$$
\left(1-\sigma_1\right)\tau_1'^-(x)-\left(1+\sigma_1\right)\nu_1^-(x)=A_1(x),\,\, 0<x<1, \eqno (21)
$$
$$
\left(1+\sigma_2\right)\tau_2'^-(t)+\left(\sigma_2-1\right)\nu_2^-(t)=A_2(t),\,\, 0<t<1, \eqno (22)
$$
$$
\left(1-\sigma_3\right)\tau_3'^+(t)+\left(1+\sigma_3\right)\nu_3^+(t)=A_3(t),\,\, 0<t<1, \eqno (23)
$$
where
$$
\begin{array}{l}
A_1(x)=2\int\limits_{\rho(x)}^x f_1(\xi_1,x)d\xi_1+2\sigma_1\int\limits_x^{v(x)} f_1(x,\eta_1)d\eta_1,\\
A_2(t)=2\int\limits_{\rho(-t)}^t f_1(\xi_1,t)d\xi_1-2\sigma_2\int\limits_t^{v(t)} f_1(t,\eta_1)d\eta_1,\\
A_3(t)=-2\int\limits_{t}^{1-v(1+t)} f_1(t,\eta_1)d\eta_1-2\sigma_3\int\limits_{\rho(1-t)-1}^{t} f_1(\xi_1, t)d\xi_1.\\
\end{array}
\eqno (24)
$$
\textbf{Remark 1.} If $|\sigma_i|=1$, then from (17)-(19) we directly can find $\tau_i'^\pm(s)$ or $\nu_i^\pm(s) (i=\overline{1,3})$. In this case, considered problem can be divided into 4 problems, which can be solved independently.

Further, we consider the case, when  $|\sigma_i|\neq 1$. Other particular cases, for instance, $|\sigma_1|=1,\, |\sigma_j|\neq 1 (j=2,3)$ will be studied similarly, but with simpler evaluations.

\section*{5. The uniqueness of the solution}

\textbf{Remark 2.} In order to prove the uniqueness of the solution of considered problem, we suppose that it has two $u_1(x,t), u_2(x,t)$ solutions. Designating difference of them as $u(x,t)=u_1(x,t)-u_2(x,t)$, we obtain corresponding homogeneous problem. Further, we prove that this homogeneous problem has only trivial solution, from which one can conclude that original problem has unique solution.

We multiply equation $u_{xx}-_CD_{0t}^\lambda u=0$ to the function $u(x,t)$ and integrate along the domain $\Omega_0$. After the usage of the Green's formulas and considering designations (20), we obtain
$$
\int\limits_0^1\tau_2^+(t)\nu_2^+(t)dt -\int\limits_0^1\tau_3^-(t)\nu_3^-(t)dt+\iint\limits_{\Omega_0} \left[u\cdot _CD_{0t}^\lambda u+u_x^2(x,t)\right]dx dt=0. \eqno (25)
$$
From (15) and (22), after some evaluations, we get
$$
\begin{array}{l}
\nu_2^+(t)=\frac{1}{1-\sigma_2}\left\{\alpha_2(1+\sigma_2)\tau_2'^+(t)+\beta_2\int\limits_t^1(1+\sigma_2)\tau_2'^+(s)P_2(t,s)ds-\right.\\
\left.-\alpha_2A_2(t)-\beta_2\int\limits_t^1A_2(s)P_2(t,s)ds\right\}.
\end{array}
\eqno (26)
$$
Similarly from (16) and (23) we have
$$
\begin{array}{l}
\nu_3^-(t)=\frac{1}{1+\sigma_3}\left\{\alpha_3(\sigma_3-1)\tau_3'^-(t)+\beta_3\int\limits_t^1(\sigma_3-1)\tau_3'^-(s)P_3(t,s)ds+\right.\\
\left.+\alpha_3A_3(t)+\beta_3\int\limits_t^1A_3(s)P_3(t,s)ds\right\}.
\end{array}
\eqno (27)
$$
Let us first investigate the sign of the integral
$$
I_1=\int\limits_0^1 \tau_2^+(t)\nu_2^+(t)dt. \eqno (28)
$$
Considering (26) at $f_1\equiv 0$ and $\tau_2(0)\equiv 0$, from (28) we get
$$
\begin{array}{l}
I_1=\frac{1+\sigma_2}{1-\sigma_2}\left\{\frac{\alpha_2}{2}\left(\tau_2^+(1)\right)^2+\beta_2\tau_2^+(1)\int\limits_0^1 \tau_2^+(t)P_2(t,1)dt-\right.\\
\left.-\beta_2\int\limits_0^1\left(\tau_2^+(t)\right)^2P_2(t,t)dt-\beta_2\int\limits_0^1\tau_2^+(t)dt\int\limits_t^1\tau_2^+(s)\frac{\partial P_2(t,s)}{\partial s}ds\right\}.
\end{array}
$$
Assuming
$$
P_2(t,s)=P_{2,1}(t)\cdot P_{2,2}(s),\,\,P_{2,2}(1)=0,
$$
from the last equality we find
$$
\begin{array}{l}
I_1=\frac{1+\sigma_2}{1-\sigma_2}\left\{-\frac{\alpha_2}{2}\left(\tau_2^+(1)\right)^2-\beta_2\int\limits_0^1\left(\tau_2^+(t)\right)^2P_{2,1}(t)P_{2,2}(t)dt\right\}-\\
-\frac{1+\sigma_2}{1-\sigma_2}\left\{\beta_2\int\limits_0^1\tau_2^+(t)dt\int\limits_t^1\tau_2^+(s)P_{2,1}(t) P_{2,2}'(s)ds\right\}=I_{1,1}-I_{1,2}.
\end{array}
$$
If
$$
\frac{\alpha_2(1+\sigma_2)}{1-\sigma_2}\geq 0,\,\,\frac{\beta_2(1+\sigma_2)}{1-\sigma_2} P_{2,1}(t) P_{2,2}(t) \leq 0, \eqno (29)
$$
then $I_{1,1}\geq 0$. Since
$$
\frac{d}{dt}\left(\int\limits_t^1\tau_2^+(s)P_{2,2}'(s)ds\right)^2=-2\tau_2^+(t)P_{2,2}'(t)\int\limits_t^1\tau_2^+(s)P_{2,2}'(s)ds,
$$
then after the integration by parts, we have
$$
\begin{array}{l}
I_{1,2}=\frac{\beta_2(1+\sigma_2)}{2(1-\sigma_2)}\left\{\frac{P_{2,1}(0)}{P_{2,2}'(0)}\left(\int\limits_0^1\tau_2^+(s)P_{2,2}'(s)ds\right)^2+\right.\\
\left.+\int\limits_0^1\left(\int\limits_0^1
\tau_2^+(s)P_{2,2}'(s)ds\right)^2\left(\frac{P_{2,1}(t)}{P_{2,2}'(t)}\right)'dt\right\}.
\end{array}
$$
If
$$
\frac{\beta_2(1+\sigma_2)}{2(1-\sigma_2)} \cdot \frac{P_{2,1}(0)}{P_{2,2}'(0)}\leq 0,\,\,\frac{\beta_2(1+\sigma_2)}{2(1-\sigma_2)} \cdot \left(\frac{P_{2,1}(t)}{P_{2,2}'(t)}\right)'\leq 0,\eqno (30)
$$
then $I_{1,2}\leq 0$. Hence, it follows that $I_1\geq 0$.

Similarly, using (27) under the certain restrictions on given functions and parameters we can prove the following inequality:
$$
I_2=\int\limits_0^1\tau_3^-(t)\nu_3^-(t)dt\leq 0.\eqno (31)
$$

Now we prove that $\tau_1^+(x)\equiv 0$ at $f(x,t)\equiv 0$. For this aim in the equation $u_{xx}-_CD_{0t}^\lambda u=0$ we pass to the limit as $t\rightarrow +0$ and considering designation (20) we have
$$
\tau_1''^+(x)-\Gamma(\lambda)\nu_1^+(x)=0. \eqno (32)
$$
We multiply equation (32) to the function $\tau_1^+(x)$ and integrate from $0$ to $1$:
$$
\int\limits_0^1 \tau_1''^+(x)\tau_1^+(x)dx-\Gamma(\lambda)\int\limits_0^1\tau_1^+(x)\nu_1^+(x)dx=0.
$$
Let us consider an integral $I_3=\int\limits_0^1\tau_1^+(x)\nu_1^+(x)dx$.
Considering (14), (20) and (21) we find that
$$
\begin{array}{l}
\nu_1^+(x)=\frac{1}{1+\sigma_1}\left\{(1-\sigma_1)\left[\alpha_1\tau_1'^+(x)+\beta_1\int\limits_x^1\tau_1'^+(s)P_1(x,s)ds\right]-\right.\\
\left.-\left[\alpha_1A_1(x)+\beta_1\int\limits_x^1A_1(s)P_1(x,s)ds\right]\right\}.\\
\end{array}
$$
At $f(x,t)\equiv 0$ we substitute obtained representation into $I_3$ and bearing condition (13) in mind, after some evaluations we get
$$
I_3=\frac{\beta_1(\sigma_1-1)}{\sigma_1+1}\left[\int\limits_0^1\left(\tau_1^+(x)\right)^2P_1(x,x)dx+\int\limits_0^1\tau_1^+(x)dx
\int\limits_x^1\tau_1^+(s)\frac{\partial P_1(x,s)}{\partial s}ds\right].\eqno (33)
$$
Let $P_1(x,s)=P_{1,1}(x)\cdot P_{1,2}(s)$. Then similarly doing the same steps as in $I_{1,2}$ from (33) we obtain
$$
\begin{array}{l}
I_3=\frac{\beta_1(\sigma_1-1)}{2(\sigma_1+1)}\left[2\int\limits_0^1\left(\tau_1^+(x)\right)^2 P_{1,1}(x)\cdot P_{1,2}(x)dx+\right.\\
\left.+\frac{P_{1,1}(0)}{P_{1,2}'(0)}\left(\int\limits_0^1\tau_1^+(s)P_{1,2}'(s)ds\right)^2+\int\limits_0^1\left(\frac{P_{1,1}(x)}{P_{1,2}'(x)}\right)'
\left(\int\limits_x^1\tau_1^+(s)P_{1,2}'(s)ds\right)^2dx\right]\\
\end{array}
$$
If
$$
\frac{\beta_1(\sigma_1-1)}{2(\sigma_1+1)}\geq 0,\,P_{1,1}(x)\cdot P_{1,2}(x)\geq 0,\, \left(\frac{P_{1,1}(x)}{P_{1,2}'(x)}\right)'\geq 0,\eqno (34)
$$
then $I_3\geq 0$.

From the other hand, if we consider (32), then
$$
I_3=\int\limits_0^1\tau_1^+(x)\cdot \frac{1}{\Gamma(\lambda)}\tau_1''^+(x)dx=-\frac{1}{2\Gamma(\lambda)}\int\limits_0^1\left(\tau_1^+(x)\right)^2dx.
$$
Since $\Gamma(\lambda)>0$ at $0<\lambda<1$, then $I_3\leq 0$. Therefore, we get that $I_3\equiv 0$, from which it follows that $\tau_1^+(x)\equiv 0$.

Based on result by A.M.Nakhushev [18], we have
$$
\iint\limits_{\Omega_0}u\cdot _CD_{0t}^\lambda u dx dt\geq 0.\eqno (35)
$$
Finally, considering $I_1\geq 0$ and (31), (35), from (25) we can conclude that $\tau_2^+(t)\equiv 0$, $\tau_3^-(t)\equiv 0$. Due to the solution (3) of FBP, we get that $u(x,t)\equiv 0$ in $\Omega_0$. Further, considering that $u(x,t)\in C(\overline{\Omega})$ we deduce that  $u(x,t)\equiv 0$ in $\overline{\Omega}$.

The uniqueness of the solution of the problem S is proved. Now we formulate our result as a theorem.

\textbf{Theorem 1.} Let
$$
P_i(t,s)=P_{i,1}(t)\cdot P_{i,2}(s), P_{j,2}(1)=0, i=\overline{1,3},\, j=2,3
$$
and conditions (29), (30), (34), and
$$
\frac{\alpha_3(\sigma_3-1)}{\sigma_3+1}\leq 0,\, \frac{\beta_3(\sigma_3-1)}{\sigma_3+1} P_{3,1}(t)\cdot P_{3,2}(t)\geq 0,
$$
$$
\frac{\beta_3(\sigma_3-1)}{2(\sigma_3+1)}\frac{P_{3,1}(0)}{P_{3,2}'(0)}\geq 0,\, \frac{\beta_3(\sigma_3-1)}{2(\sigma_3+1)}\left(\frac{P_{3,1}(t)}{P_{3,2}'(t)}\right)'\geq 0
$$
be hold. If there exist solution of the problem S, then it is unique.

\section*{6. The existence of the solution}

From (21) and (32) we find
$$
\tau_1''^+(x)-C\tau_1'^+(x)=F_1(x),\eqno (36)
$$
where $C=\frac{\alpha_1\Gamma(\lambda)(\sigma_1-1)}{\sigma_1+1}$,
$$
F_1(x)=-C\beta_1\int\limits_x^1\tau_1'^+(s)P_1(x,s)ds-\frac{\alpha_1\Gamma(\lambda)}{1+\sigma_1}A_1(x)-
\frac{\beta_1\Gamma(\lambda)}{1+\sigma_1}\int\limits_x^1A_1(s)P_1(x,s)ds.\eqno (37)
$$
Solution of the equation (36) with conditions $\tau_1^+(0)=\tau_1^+(1)=0$ has a form
$$
\tau_1^+(x)=\int\limits_0^1G_0(x,\xi)F_1(\xi)d\xi,
$$
where
$$
G_0(x,\xi)=\frac{1}{C\left[e^{Cx}-e^{C(x-1)}\right]}
\left\{
\begin{array}{l}
\left(1-e^{C\xi}\right)\left(1-e^{C(x-1)}\right),\,\,0\leq \xi\leq x,\\
\left(1-e^{C(\xi-1)}\right)\left(1-e^{Cx}\right),\,\,x\leq \xi\leq 1.\\
\end{array}
\right.
$$
Considering (36) and (37) we obtain integral equation
$$
\tau_1^+(x)-\int\limits_0^1 \tau_1^+(\eta)K(x,\eta)d\eta=\overline{F_1}(x),\eqno (38)
$$
where
$$
K(x,\eta)=C\beta_1\left[G_0(x,\eta)P_1(\eta,\eta)+\int\limits_\eta^1G_0(x,\xi)\frac{\partial P_1(\xi,\eta)}{\partial \eta}d\xi\right],
$$
$$
\overline{F_1}(x)=-\frac{\Gamma(\lambda)}{1+\sigma_1}\int\limits_0^1G_(x,\eta)\left[\alpha_1A_1(x)-\beta_1\int\limits_\eta^1A_1(s)P_1(\eta,s)ds\right]d\eta.
$$
Supposing $f(\cdot,\cdot),\,P_1(\cdot,\cdot)\in C\left([0,1]\times [0,1]\right)$ and based on general theory of Fredholm integral equations, solution of the equation (38) we can write via resolvent-kernel:
$$
\tau_1^+(x)=\overline{F_1}(x)-\int\limits_0^1 \overline{F_1}(\eta)R(x,\eta)d\eta,\eqno (39)
$$
where $R(x,\eta)$ is resolvent of the kernel $K(x,\eta)$.

Further, since $\tau_i^+(x)=\tau_i^-(x),\,\,\tau_i'^+(x)=\tau_i'^-(x)\,\,(i=\overline{1,3})$, we omit signs in upper indexes.
Based on (3) we write solution of FBP as follows
$$
\begin{array}{l}
u(x,t)=\int\limits_0^t G_\xi(x,t,0,\eta)\tau_2(\eta)d\eta-\int\limits_0^t G_\xi(x,t,1,\eta)\tau_3(\eta)d\eta+\\
+\int\limits_0^1 \overline{G}(x-\xi,t)\tau_1(\xi)d\xi-\int\limits_0^1 \int\limits_0^t G(x,t,\xi,\eta) f(\xi, \eta)d\xi d\eta.\\
\end{array}
\eqno (40)
$$
Respectively, Green's function of FBP has a form
$$
G(x,t,\xi,\eta)=\frac{(t-\eta)^{\beta-1}}{2}\sum_{n=-\infty}^{+\infty} \left[e_{1,\beta}^{1,\beta} \left(-\frac{|x-\xi+2n|}{(t-\eta)^\beta}\right)
-e_{1,\beta}^{1,\beta} \left(-\frac{|x+\xi+2n|}{(t-\eta)^\beta}\right)\right].\eqno (41)
$$

We differentiate (41) once on $\xi$ and once on $x$, then passing to the limits as $x\rightarrow +0$ and $x\rightarrow 1-0$, after some evaluations we find
$$
G_{\xi x}(+0,t,0,\eta)=\frac{\partial}{\partial\eta}\left(\sum\limits_{n=-\infty}^{+\infty}\frac{1}{(t-\eta)^\beta}e_{1,1-\beta}^{1,\beta} \left(-\frac{|2n|}{(t-\eta)^\beta}\right)\right),\eqno (42)
$$
$$
G_{\xi x}(+0,t,1,\eta)=\frac{\partial}{\partial\eta}\left(\sum\limits_{n=-\infty}^{+\infty}\frac{1}{(t-\eta)^\beta}e_{1,1-\beta}^{1,\beta} \left(-\frac{|2n+1|}{(t-\eta)^\beta}\right)\right),\eqno (43)
$$
$$
G_{\xi x}(1-0,t,0,\eta)=\frac{\partial}{\partial\eta}\left(\sum\limits_{n=-\infty}^{+\infty}\frac{1}{(t-\eta)^\beta}e_{1,1-\beta}^{1,\beta} \left(-\frac{|2n+1|}{(t-\eta)^\beta}\right)\right),\eqno (44)
$$
$$
G_{\xi x}(1-0,t,1,\eta)=\frac{\partial}{\partial\eta}\left(\sum\limits_{n=-\infty}^{+\infty}\frac{1}{(t-\eta)^\beta}e_{1,1-\beta}^{1,\beta} \left(-\frac{|2n|}{(t-\eta)^\beta}\right)\right).\eqno (45)
$$
We differentiate (40) with respect to $x$ and get
$$
\begin{array}{l}
u_x(x,t)=\int\limits_0^t G_{\xi x}(x,t,0,\eta)\tau_2(\eta)d\eta-\int\limits_0^t G_{\xi x}(x,t,1,\eta)\tau_3(\eta)d\eta+\\
+\int\limits_0^1 \overline{G}_x(x-\xi,t)\tau_1(\xi)d\xi-\int\limits_0^1 \int\limits_0^t G_x(x,t,\xi,\eta) f(\xi, \eta)d\xi d\eta.\\
\end{array}
\eqno (46)
$$
At $x\rightarrow+0$ considering (42), (43) from (46) we deduce
$$
\begin{array}{l}
u_x(+0,t)=\nu_2^+(t)=\int\limits_0^t \tau_2(\eta)\frac{\partial}{\partial\eta}\left(\sum\limits_{n=-\infty}^{+\infty}\frac{1}{(t-\eta)^\beta}e_{1,1-\beta}^{1,\beta} \left(-\frac{|2n|}{(t-\eta)^\beta}\right)\right)d\eta-\\
-\int\limits_0^t \tau_3(\eta)\frac{\partial}{\partial\eta}\left(\sum\limits_{n=-\infty}^{+\infty}\frac{1}{(t-\eta)^\beta}e_{1,1-\beta}^{1,\beta} \left(-\frac{|2n+1|}{(t-\eta)^\beta}\right)\right)d\eta+\\
+\int\limits_0^1 \overline{G}_x(-\xi,t)\tau_1(\xi)d\xi-\int\limits_0^1 \int\limits_0^t G_x(0,t,\xi,\eta) f(\xi, \eta)d\xi d\eta.\\
\end{array}
\eqno (47)
$$
Similarly, as $x\rightarrow 1-0$, considering (44), (45) from (46) we have
$$
\begin{array}{l}
u_x(1-0,t)=\nu_3^-(t)=\int\limits_0^t \tau_2(\eta)\frac{\partial}{\partial\eta}\left(\sum\limits_{n=-\infty}^{+\infty}\frac{1}{(t-\eta)^\beta}e_{1,1-\beta}^{1,\beta} \left(-\frac{|2n+1|}{(t-\eta)^\beta}\right)\right)d\eta-\\
-\int\limits_0^t \tau_3(\eta)\frac{\partial}{\partial\eta}\left(\sum\limits_{n=-\infty}^{+\infty}\frac{1}{(t-\eta)^\beta}e_{1,1-\beta}^{1,\beta} \left(-\frac{|2n|}{(t-\eta)^\beta}\right)\right)d\eta+\\
+\int\limits_0^1 \overline{G}_x(1-\xi,t)\tau_1(\xi)d\xi-\int\limits_0^1 \int\limits_0^t G_x(1,t,\xi,\eta) f(\xi, \eta)d\xi d\eta.\\
\end{array}
\eqno (48)
$$
In (47), (48) we use formula of integration by parts and considering $\tau_2(0)=\tau_3(0)=0$, which follows from (13), and obtain
$$
\begin{array}{l}
\nu_2^+(t)=-\int\limits_0^t \tau_2'(\eta)K_1(t,\eta)d\eta+\int\limits_0^t \tau_3'(\eta)K_2(t,\eta)d\eta+\\
+\int\limits_0^1 \overline{G}_x(-\xi,t)\tau_1(\xi)d\xi-\int\limits_0^1 \int\limits_0^t G_x(0,t,\xi,\eta) f(\xi, \eta)d\xi d\eta,\\
\end{array}
\eqno (49)
$$
$$
\begin{array}{l}
\nu_3^-(t)=-\int\limits_0^t \tau_2'(\eta)K_2(t,\eta)d\eta+\int\limits_0^t \tau_3'(\eta)K_1(t,\eta)d\eta+\\
+\int\limits_0^1 \overline{G}_x(1-\xi,t)\tau_1(\xi)d\xi-\int\limits_0^1 \int\limits_0^t G_x(1,t,\xi,\eta) f(\xi, \eta)d\xi d\eta,\\
\end{array}
\eqno (50)
$$
where
$$
\begin{array}{l}
K_1(t,\eta)=\sum\limits_{n=-\infty}^{+\infty}\frac{1}{(t-\eta)^\beta}e_{1,1-\beta}^{1,\beta} \left(-\frac{|2n|}{(t-\eta)^\beta}\right),\\
K_2(t,\eta)=\sum\limits_{n=-\infty}^{+\infty}\frac{1}{(t-\eta)^\beta}e_{1,1-\beta}^{1,\beta} \left(-\frac{|2n+1|}{(t-\eta)^\beta}\right).\\
\end{array}
\eqno (51)
$$
From (22), (23) we find functions $\nu_2^+(t),\,\nu_3^-(t)$ and substituting them into transmitting conditions (15), (16) we find
$$
\begin{array}{l}
\nu_2^+(t)=\frac{\alpha_2(1+\sigma_2)}{1-\sigma_2}\tau_2'(t)+\frac{\beta_2(1+\sigma_2)}{1-\sigma_2}\int\limits_t^1 \tau_2'(s)P_2(t,s)ds+\frac{\alpha_2A_2(t)}{\sigma_2-1}+\\
+\frac{\beta_2}{\sigma_2-1}\int\limits_t^1A_2(s)P_2(t,s)ds,\\
\end{array}
\eqno (52)
$$
$$
\begin{array}{l}
\nu_3^-(t)=\frac{\alpha_3(\sigma_3-1)}{\sigma_3+1}\tau_3'(t)+\frac{\beta_3(\sigma_3-1)}{\sigma_3+1}\int\limits_t^1 \tau_3'(s)P_3(t,s)ds+\frac{\alpha_3A_3(t)}{1+\sigma_3}+\\
+\frac{\beta_3}{1+\sigma_3}\int\limits_t^1A_3(s)P_3(t,s)ds.\\
\end{array}
\eqno (53)
$$
(52), (53) we substitute into (49), (50), respectively:
$$
\begin{array}{l}
\frac{\alpha_2(1+\sigma_2)}{1-\sigma_2}\tau_2'(t)+ \int\limits_0^t \tau_2'(\eta)K_1(t,\eta)d\eta+\int\limits_t^1 \tau_2'(\eta)\frac{\beta_2(1+\sigma_2)}{1-\sigma_2}P_2(t,\eta)d\eta=\\
=\int\limits_0^t \tau_3'(\eta)K_2(t,\eta)d\eta+\frac{\alpha_2A_2(t)}{1-\sigma_2}+\frac{\beta_2}{1-\sigma_2}\int\limits_t^1A_2(\eta)P_2(t,\eta)d\eta+\\
+\int\limits_0^1\overline{G}_x(-\xi,t)\tau_1(\xi)d\xi-\int\limits_0^1\int\limits_0^tG_x(0,t,\xi,\eta)f(\xi,\eta)d\xi d\eta,\\
\end{array}
\eqno (54)
$$
$$
\begin{array}{l}
\frac{\alpha_3(\sigma_3-1)}{\sigma_3+1}\tau_3'(t)- \int\limits_0^t \tau_3'(\eta)K_1(t,\eta)d\eta-\int\limits_t^1 \tau_3'(\eta)\frac{\beta_3(1-\sigma_3)}{1+\sigma_3}P_3(t,\eta)d\eta=\\
=-\int\limits_0^t \tau_2'(\eta)K_2(t,\eta)d\eta-\frac{\alpha_3A_3(t)}{1+\sigma_3}-\frac{\beta_3}{1+\sigma_3}\int\limits_t^1A_3(\eta)P_3(t,\eta)d\eta+\\
+\int\limits_0^1\overline{G}_x(1-\xi,t)\tau_1(\xi)d\xi-\int\limits_0^1\int\limits_0^tG_x(1,t,\xi,\eta)f(\xi,\eta)d\xi d\eta.\\
\end{array}
\eqno (55)
$$

Let $\alpha_2\neq 0,\,\alpha_3\neq 0$. Then from (54) and (55) we find
$$
\left\{
\begin{array}{l}
\tau_2'(t)+\int\limits_0^1 \tau_2'(\eta)K_3(t,\eta)d\eta=F_2(t),\\
\tau_3'(t)+\int\limits_0^1 \tau_3'(\eta)K_4(t,\eta)d\eta=F_3(t),
\end{array}
\right.
\eqno (56)
$$
where
$$
\begin{array}{l}
K_3(t,\eta)=\frac{1-\sigma_2}{\alpha_2(1+\sigma_2)}\left\{
\begin{array}{l}
K_1(t,\eta),\,\,\,\,\,\,\,\,\,\,\,\,\,0\leq\eta\leq t,\\
\frac{\beta_2(1+\sigma_2)}{1-\sigma_2}P_2(t,\eta),\,t\leq\eta\leq 1,\\
\end{array}
\right.\\
K_4(t,\eta)=\frac{\sigma_3+1}{\alpha_3(\sigma_3-1)}\left\{
\begin{array}{l}
K_1(t,\eta),\,\,\,\,\,\,\,\,\,\,\,\,\,0\leq\eta\leq t,\\
\frac{\beta_3(1-\sigma_3)}{1+\sigma_3}P_3(t,\eta),\,t\leq\eta\leq 1,\\
\end{array}
\right.
\end{array}
\eqno (57)
$$
$$
\begin{array}{l}
F_2(t)=\frac{1-\sigma_2}{\alpha_2(1+\sigma_2)}\left\{\int\limits_0^t\tau_3'(\eta)K_2(t,\eta)d\eta+\frac{\alpha_2A_2(t)}{1-\sigma_2}+\right.\\
\left.+\int\limits_t^1\frac{\beta_2A_2(\eta)P_2(t,\eta)}{1-\sigma_2}d\eta+\int\limits_0^1\overline{G}_x(-\xi,t)\tau_1(\xi)d\xi-
\int\limits_0^1\int\limits_0^tG_x(0,t,\xi,\eta)f(\xi,\eta)d\xi d\eta\right\},
\end{array}
$$
$$
\begin{array}{l}
F_3(t)=\frac{\sigma_3+1}{\alpha_3(\sigma_3-1)}\left\{-\int\limits_0^t\tau_2'(\eta)K_2(t,\eta)d\eta-\frac{\alpha_3A_3(t)}{1+\sigma_3}-\right.\\
\left.-\int\limits_t^1\frac{\beta_3A_3(\eta)P_3(t,\eta)}{1+\sigma_3}d\eta+\int\limits_0^1\overline{G}_x(1-\xi,t)\tau_1(\xi)d\xi-
\int\limits_0^1\int\limits_0^tG_x(1,t,\xi,\eta)f(\xi,\eta)d\xi d\eta\right\},
\end{array}
$$
functions $A_i(t) (i=\overline{1,3})$ are defined by (24).

We formally represent solution of the first equation of the system (56) via resolvent-kernel:
$$
\tau_2'(t)=F_2(t)+\int\limits_0^1 F_2(\eta)R_1(t,\eta)d\eta,\eqno (58)
$$
where $R_1(t,\eta)$ is resolvent of the kernel $K_3(t,\eta)$.
Considering representation of $F_2(t)$, (58) substitute into the representation of $F_3(t)$. After some evaluations obtained representation of $F_3(t)$ we substitute into the second equation of the system (56):
$$
\tau_3'(t)+\int\limits_0^1 \tau_3'(\eta)K_5(t,\eta)d\eta=F_4(t),\eqno (59)
$$
where
$$
K_5(t,\eta)=K_4(t,\eta)-\int\limits_\eta^1\frac{(\sigma_3+1)(\sigma_2-1)}{\alpha_2\alpha_3(\sigma_2+1)(\sigma_3-1)}\overline{K}_2(t,s)K_2(s,\eta)ds,
$$
$$
\overline{K}_2(t,s)=\left\{
\begin{array}{l}
K_2(t,s)+\int\limits_0^tR_2(z,s)K_2(t,z)dz,\,0\leq s\leq t,\\
\int\limits_0^tR_1(z,s)K_2(t,z)dz,\,\,\,\,\,\,\,\,\,\,\,\,\,\,\,t\leq s\leq 1,\\
\end{array}
\right.
$$
$$
\begin{array}{l}
F_4(t)=\frac{(\sigma_3+1)(1-\sigma_2)}{\alpha_2\alpha_3(\sigma_2+1)(\sigma_3-1)}\left\{-\frac{\alpha_2}{1-\sigma_2}
\int\limits_0^1\overline{K}_2(t,\eta)d\eta\int\limits_0^1A_2(s)P_2(\eta,s)ds-\right.\\
\left.-\int\limits_0^1\overline{K}_2(t,\eta)d\eta\int\limits_0^1\overline{G}_x(-\xi,\eta)\tau_1(\xi)d\xi+
\int\limits_0^1\overline{K}_2(t,\eta)d\eta\int\limits_0^1\int\limits_0^\eta G_x(0,\eta,\xi,s) f(\xi,s)d\xi ds\right\}+\\
+\frac{\sigma_3+1}{\alpha_3(\sigma_3-1)}\left\{-\frac{\alpha_3A_3(t)}{1+\sigma_3}-\frac{\beta_3}{1+\sigma_3}\int\limits_t^1A_3(\eta)P_3(t,\eta)d\eta+\right.\\
\left.+\int\limits_0^1\overline{G}_x(1-\xi,t)\tau_1(\xi)d\xi-\int\limits_0^1\int\limits_0^t G_x(1,\eta,\xi,\eta) f(\xi,\eta)d\xi d\eta\right\}.
\end{array}
$$
Since $F_4(t)\in C[0,1]\cap C^1(0,1)$ and $|K_5(t,\eta)|\leq\frac{C}{(t-\eta)^\beta}$, due to (51), (57), we rewrite solution of (59) via resolvent-kernel:
$$
\tau_3'(t)=F_4(t)+\int\limits_0^1 F_4(\eta)R_2(t,\eta)d\eta,\eqno (60)
$$
where $R_2(t,\eta)$ is resolvent of the kernel $K_5(t,\eta)$.
Further, using transmitting conditions (15), (16) and main functional correlations (22), (23) we can find functions $\nu_j^\pm(t) (j=2,3)$. Based on solution of FBP and Cauchy problems, we can recover solution of the problem by formulas (17)-(19), (40).

Let us know consider the case, when $\alpha_j=0 (j=2,3)$.

Based on representation of $K_1(t,\eta)$ at $n=0$, from (54) we get
$$
\int\limits_0^t \frac{\tau_2'(\eta)}{(t-\eta)^\beta}d\eta=\overline{F}_2(t),\eqno (61)
$$
where
$$
\begin{matrix}
\overline{F}_2(t)=\frac{\beta_2(\sigma_2+1)}{\sigma_2-1}\int\limits_t^1\tau_2'(\eta)P_2(t,\eta)d\eta-\int\limits_0^t\tau_2'(\eta)\overline{K}_1(t,\eta)d\eta+\\
+\int\limits_0^t\tau_3'(\eta)K_2(t,\eta)d\eta+\frac{\alpha_2A_2(t)}{1-\sigma_2}+\frac{\beta_2}{1-\sigma_2}\int\limits_t^1A_2(\eta)P_2(t,\eta)d\eta+\\
+\int\limits_0^1\overline{G}_x(-\xi,t)\tau_1(\xi)d\xi-\int\limits_0^1\int\limits_0^tG_x(0,t,\xi,\eta)f(\xi,\eta)d\xi d\eta,\\
\end{matrix}
$$
$$
\overline{K}_1(t,\eta)=\sum\limits_{\begin{smallmatrix} n=-\infty\\n\neq 0 \end{smallmatrix}}^\infty\frac{1}{(t-\eta)^\beta}e_{1,\beta}^{1,1-\beta}\left(-\frac{|2n|}{(t-\eta)^\beta}\right).
$$
(61) is a generalized Abel's integral equation, solution of which we can rewrite as follows [19]
$$
\tau_2'(t)=\frac{\sin\beta\pi}{\pi}\left[\frac{\overline{F}(0)}{t^{1-\beta}}+\int\limits_0^t \frac{\overline{F}_2'(z)dz}{(t-z)^{1-\beta}}\right].
$$
Considering representation of $\overline{F}_2(t)$ after some evaluations we deduce
$$
\tau_2'(t)-\int\limits_0^1\tau_2'(\eta)K_6(t,\eta)d\eta=F_5(t),\eqno (62)
$$
where
$$
K_6(t,\eta)=\frac{\beta_2(\sigma_2+1)\sin\beta\pi}{\pi(\sigma_2-1)t^{1-\beta}}
\left\{
\begin{matrix}
P_2(0,\eta)-\frac{P_2(\eta,\eta)}{(t-\eta)^{1-\beta}}
+\int\limits_\eta^t\frac{\partial P_2(z,\eta)}{\partial z}\frac{dz}{(t-z)^{1-\beta}}-\\
-\frac{\sigma_2-1}{\beta_2(\sigma_2+1)}\int\limits_\eta^t \frac{\partial \overline{K}_1(z,\eta)}{\partial z}\frac{dz}{(t-z)^{1-\beta}},\,\,0\leq \eta\leq t,\\
P_2(0,\eta)+\int\limits_0^t\frac{\partial P_2(z,\eta)}{\partial z}\frac{dz}{(t-z)^{1-\beta}},\,\,t\leq \eta\leq 1,\\
\end{matrix}
\right.
$$
$$
\begin{array}{l}
F_5(t)=\frac{\sin\beta\pi}{\pi}\left\{\int\limits_0^1\tau_3'(\eta)d\eta\int\limits_0^t\frac{\partial K_2(z,\eta)}{\partial z}\frac{dz}{(t-z)^{1-\beta}}+\right.\\
+\frac{1}{t^{1-\beta}}\left[\frac{\alpha_2A_2(0)}{1-\sigma_2}+\frac{\beta_2}{1-\sigma_2}\int\limits_0^1A_2(\eta)P_2(0,\eta)d\eta+\int\limits_0^1\overline{G}_x(-\xi,0)\tau_1(\xi)d\xi\right]+\\
+\frac{\alpha_2}{1-\sigma_2}\int\limits_0^t\frac{A_2'(z)dz}{(t-z)^{1-\beta}}-\frac{\beta_2}{1-\sigma_2}\int\limits_0^t\frac{A_2(z)P_2(z,z)}{(t-z)^{1-\beta}}dz+
\frac{\beta_2}{1-\sigma_2}\int\limits_0^t\frac{dz}{(t-z)^{1-\beta}}\int\limits_z^1A_2(\eta)\frac{\partial P_2(z,\eta)}{\partial z}d\eta+\\
+\int\limits_0^t\frac{dz}{(t-z)^{1-\beta}}\int\limits_0^1\overline{G}_{xz}(-\xi,z)\tau_1(\xi)d\xi-\int\limits_0^t\frac{dz}{(t-z)^{1-\beta}}\int\limits_0^1G_x(0,z,\xi,z)f(\xi,z)d\xi-\\
\left.-\int\limits_0^t\frac{dz}{(t-z)^{1-\beta}}\int\limits_0^1\int\limits_0^z G_{xz}(0,z,\xi,\eta)f(\xi,\eta)d\xi d\eta\right\}.
\end{array}
$$
Based on (51) and supposing $F_5(t)$ as known, solution of (62) we write via resolvent:
$$
\tau_2'(t)=F_5(t)-\int\limits_0^1F_5(\eta)R_3(t,\eta)d\eta
$$
and substitute it into the representation of $F_3(t)$, and after some evaluations from the second equation of the system (56) we deduce
$$
\tau_3'(t)-\int\limits_0^1\tau_3'(\eta)K_7(t,\eta)d\eta=F_6(t),\eqno (63)
$$
where
$$
\begin{array}{l}
K_7(t,\eta)=K_4(t,\eta)+\frac{1+\sigma_3}{\alpha_3(1-\sigma_3)}\frac{\sin\beta\pi}{\pi}\left[\int\limits_0^tK_2(t,z)dz\int\limits_0^z\frac{\partial K_2(s,\eta)}{\partial s}\frac{ds}{(z-s)^{1-\beta}}-\right.\\
\left.-\int\limits_0^tK_2(t,z)dz\int\limits_0^1R_2(z,s)ds\int\limits_0^s \frac{\partial K_2(\delta,\eta)}{\partial \delta}\frac{d\delta}{(\eta-\delta)^{1-\beta}}\right],
\end{array}
$$
$$
\begin{array}{l}
F_6(t)=\frac{1+\sigma_3}{\alpha_3(1-\sigma_3)}\frac{\sin\beta\pi}{\pi}\left[\int\limits_0^t\overline{F}_5(\eta)K_2(t,\eta)d\eta-\int\limits_0^1\overline{F}_5(z)dz
\int\limits_0^t R_3(\eta,z)K_2(t,\eta)d\eta\right],\\
\overline{F}_5(t)=\frac{\sin\beta\pi}{\pi}\left\{\frac{1}{t^{1-\beta}}\left[\frac{\alpha_2A_2(0)}{1-\sigma_2}+\frac{\beta_2}{1-\sigma_2}\int\limits_0^1A_2(\eta)P_2(0,\eta)d\eta+\right.\right.\\
\left.+\int\limits_0^1\overline{G}_x(-\xi,0)\tau_1(\xi)d\xi\right]+\frac{\alpha_2}{1-\sigma_2}\int\limits_0^t\frac{A_2'(z)dz}{(t-z)^{1-\beta}}-\frac{\beta_2}{1-\sigma_2}
\int\limits_0^1\frac{A_2(z)P_2(z,z)}{(t-z)^{1-\beta}}dz+\\
+\frac{\beta_2}{1-\sigma_2}\int\limits_0^t\frac{dz}{(t-z)^{1-\beta}}\int\limits_z^1A_2(\eta)\frac{\partial P_2(z,\eta)}{\partial z}d\eta+\int\limits_0^t\frac{dz}{(t-z)^{1-\beta}}\int\limits_0^1\overline{G}_{xz}(-\xi,z)\tau_1(\xi)d\xi-\\
\left.-\int\limits_0^t\frac{dz}{(t-z)^{1-\beta}}\int\limits_0^1G_x(0,z,\xi,z)f(\xi,z)d\xi-\int\limits_0^t\frac{dz}{(t-z)^{1-\beta}}\int\limits_0^1\int\limits_0^zG_{xz}(0,z,\xi,\eta)f(\xi,\eta)d\xi d\eta\right\}.
\end{array}
$$
Imposing certain conditions to given functions we can state that $|K_7(t,\eta)|\leq\frac{C}{|t-\eta|^{1-\beta}}$ and $F_6(t)\in C^1(0,1)$. Further, as in the previous case, we can find solution of (63) via resolvent:
$$
\tau_3'(t)=F_6(t)-\int\limits_0^1F_6(\eta)R_4(t,\eta)d\eta,
$$
where $R_4(t,\eta)$ is resolvent of the kernel $K_7(t,\eta)$.

Since we found functions $\tau_i(\cdot) (i=\overline{1,3})$, functions $\nu_i^\pm(\cdot)$ can be found by formulas (39), (58), (60), (21)-(23) and transmitting conditions (14)-(16).

Solution of the problem S in $\Omega_0$ we will recover using the solution of FBP by the formula (40), in $\Omega_i (i=\overline{1,3})$ by the formulas (17)-(19).

We proved the following existence theorem:

{\bf Theorem 2.} If
$$
f(x,t),\,P_i(x,t)\in C\left([0,1]\times [0,1]\right)\cap C^1\left((0,1)\times (0,1)\right),
$$
then there exists solution of the problem S.

\section*{7. Conclusion.}
Below we highlight distinctive sides of considered problem:
\begin{itemize}
	\item Parabolic-hyperbolic type equation contains Caputo fractional operator with respect to time variable in a parabolic part. As a consequence of this we had to deal with properties of Wright type function in order to simplify main functional correlations;
	\item Considered mixed domain bounded by curves, on which non-local conditions were given. Under appropriate assumptions to the curves, we obtain simpler form of main functional correlations, which make further evaluations possible. 
	\item Transmitting conditions on the type-changing lines have integral form, which lead to the separation of investigaion of the problem to different cases. Presicely, in one case directly, in other case, using solution of general Abel's integral equation we reduce considered problem to the system of Fredholm integral equations.
\end{itemize}

For the proof of the uniqueness result (see Theorem 1) we mainly use energy integrals with appropriate modifications. The proof of the existence theorem (see Theorem 2), we realize by reducing the problem to the system of Fredholm integral equations.

Obtained results will give a possibility to study spectral properties of such problems. On this direction we refer works [20-21].

\section*{Acknowledgement}
This work partially supported by the Grant No 3293/GF4 of the Ministry of education and science of the Republic of Kazakhstan.

\end{document}